\documentclass[a4paper, reqno,12pt]{amsart}
\usepackage{tikz, tikz-3dplot}
\usepackage{amsmath, amsthm, amssymb}

\theoremstyle{plain}
\newtheorem{te}{Theorem}[section]
\newtheorem{lem}[te]{Lemma}

\newtheorem{pr}[te]{Proposition}

\newtheorem{qu}[te]{Question}
\newtheorem{con}[te]{Conjecture}
\theoremstyle{remark}
\newtheorem{re}[te]{Remark}
\newtheorem*{ack*}{Acknowledgment}

\def\x{{\textbf{x}}}

\def\n{{\bf n}}

\def\0{{\bf 0}}

\def\T{{\mathbb T}}
\def\R{{\mathbb R}}

\def\N{{\mathbb N}}
\def\C{{\mathbb C}}
\def\S{{\mathbb S}}
\def\Z{{\mathbb Z}}

\def\nint{\mathop{\diagup\kern-13.0pt\int}}

\def\les{{\;\lessapprox}\;}

\begin{document}

\title[$L^2$ restriction ]{On restriction of exponential sums to hypersurfaces with zero curvature}

\author{Ciprian Demeter}
\address{Department of Mathematics, Indiana University, 831 East 3rd St., Bloomington IN 47405}
\email{demeterc@indiana.edu}

\thanks{The author is  partially supported by the NSF grant DMS-2055156}

\begin{abstract}
We prove essentially sharp bounds for the $L^p$ restriction of weighted Gauss sums to monomial curves. Getting the $L^2$ upper bound combines the $TT^*$ method for matrices with the first and second derivative test for exponential sums. The matching lower bound follows via constructive interference on short blocks of integers, near the critical point of the phase function. This method is used to make the broader point that restriction to hypersurfaces is really sensitive to curvature. Our results here complement those in \cite{DL}.  	
\end{abstract}

\maketitle

\section{introduction}

The paper \cite{DL} is concerned with proving essentially scale-independent (by that we mean estimates with at most $N^\epsilon$ losses) $L^p$ bounds for a variety of weighted exponential sums restricted to hypersurfaces. The key assumption in \cite{DL} is that hypersurfaces have nonzero Gaussian curvature. This property was exploited in our arguments by means of the decay of the Fourier transform of the corresponding surface measures.

The simplest instance is provided by the weighted Gauss sums. We write $e(z)=e^{2\pi iz}$.
\begin{te}[\cite{DL}]
\label{23}	
Assume $\phi:[-1,1]\to\R$ is $C^2$ and satisfies $\inf_{x\in[-1,1]}|\phi''(x)|>0$. Then for each $a_n\in \C$ we have
$$\|\sum_{n=-N}^Na_ne(nx+n^2\phi(x))\|_{L^2(-1,1)}\lesssim\|a_n\|_{2}.$$
\end{te}

For all the examples discussed in \cite{DL}, it was not clear whether nonzero curvature was an artifact of the proof or a necessary condition. Here, we prove the latter is the case. For the example of weighted Gauss sums restricted to monomial curves, we can in fact be very precise.

\begin{te}
\label{20}Let $k\ge 3$. Let $\mathcal{B}_{N,k}$ the smallest constant such that
the inequality
$$\|\sum_{n=-N}^Na_ne(nx-n^2x^k)\|_{L^2(-1,1)}\le \mathcal{B}_{N,k}\|a_n\|_{2}$$
holds for arbitrary $a_n\in\C$. Then $N^{\frac{k-2}{6(k-1)}}\lesssim \mathcal{B}_{N,k}\les N^{\frac{k-2}{6(k-1)}}$.
\end{te}
The notation $\les$ stands for logarithmic losses, $(\log N)^{O(1)}$ in the scale parameter $N$.
\smallskip

The main idea for getting lower bounds is described in the next section. We construct coefficients $a_n$ that realize constructive interference for our weighted Gauss sums on a ``large" set of values for $x$. Standard computations show that the constant coefficient case $a_n=1$ is well behaved, in the sense that
$$\|\sum_{n=-N}^Ne(nx-n^2x^k)\|_{L^2(-1,1)}\lesssim_\epsilon N^{\frac12+\epsilon}.$$
Our optimal construction will involve quadratically modulated coefficients. They arise while doing Taylor expansions near the critical point of the phase function. Similar scale-dependent lower bounds in higher dimensions are explored in Section \ref{s4}.

The matching upper bound in  Theorem  \ref{20} is obtained in Section \ref{s3}. The $TT^*$ method is combined with stationary phase analysis and derivative  tests for exponential sums.

The numerology in Theorem \ref{20} is interesting in a few ways. The sharp bound at the endpoint $L^2$ is sensitive to the exponent $k\ge 2$, but the sharp bound at the critical exponent $p=4$ is not. See Theorem \ref{44}. Also, both the lower and  upper bounds in $L^2$ are more difficult to obtain than the ones in $L^4$. This is unlike the typical examples in the literature.
\smallskip

The result in Theorem \ref{20} can also be related to the following (essentially equivalent reformulation of the) periodic  Schr\"odinger maximal operator conjecture. For each measurable function $\phi:[0,1]\to [0,1]$, it is expected that the following estimate holds
$$\|\sum_{n=-N}^Na_ne(nx+n^2\phi(x))\|_{L^2(0,1)}\lesssim_\epsilon N^{\frac14+\epsilon} \|a_n\|_2.$$
See \cite{Ba} for recent progress and the history of this conjecture.

\begin{ack*}
	The author would like to thank Alex Barron for helpful discussions on $L^4$, and the referee for a very careful reading of the manuscript.
\end{ack*}

\section{Lower bounds}
\label{s2}
Throughout the rest of the paper we will write either $A=o(B)$ or $A\ll B$ to denote the fact that $A\le cB$ for some small enough absolute constant $c>0$, that depends on $k$. The notation $A\lesssim B$ or $A=O(B)$ will refer to the scenario when $A\le CB$ for some potentially large, but still universal, constant $C$.

\begin{te}
\label{21}	
Assume $\phi\in C^{k+1}(-1,1)$ for some $k\ge 3$. Assume there is $x_0\in(-1,1)$ such that
	$$\phi'(x_0)=\phi''(x_0)=\ldots=\phi^{(k-1)}(x_0)=0$$
	$$\phi^{(k)}(x_0)\not=0.$$
	Then for each $N$  there exists $a_n\in \C$ such that for each $p\ge 2$
$$\|\sum_{n=-N}^Na_ne(nx-n^2\phi(x))\|_{L^p(-1,1)}\gtrsim N^{\frac{2k-1}{6(k-1)}-\frac{k+1}{3p(k-1)}}\|a_n\|_2, $$
with an implicit constant depending on $\phi$, but independent of $N$.
\end{te}
\begin{proof}
Due to modulation invariance, we may assume $x_0=0$,  at the expense of replacing the integration domain with a smaller interval $(-\delta,\delta)$. We will also normalize such that  $\phi(x)={x^k}+O(x^{k+1})$. Indeed, since the sum in $n$ is over a symmetric interval, it will not matter whether $\phi^{(k)}(0)$ is positive or negative. Let $x_N\sim\frac{1}{N^{\frac1{k-1}}}$ be such that $\phi'(x_N)=\frac1N$. Such a point exists since $\phi'(x)=kx^{k-1}+O(x^k)$.

 We apply Taylor's formula to $\psi(x,n)=nx-n^2\phi(x)$ near the point $(x_N,\frac{N}2)$. If $x\in I_1=[x_N-\Delta, x_N+\Delta]$ and $n\in I_2=[\frac{N}{2}-M,\frac{N}{2}+M]$ then
\begin{align*}
\psi(x,n)&=\psi(x_N,\frac{N}2)+\psi_x(x_N,\frac{N}2)(x-x_N)+\psi_n(x_N,\frac{N}2)(n-\frac{N}{2})\\&+\frac12\psi_{xx}(x_N,\frac{N}2)(x-x_N)^2+\frac12\psi_{nn}(x_N,\frac{N}2)(n-\frac{N}{2})^2\\&+\psi_{xn}(x_N,\frac{N}2)(x-x_N)(n-\frac{N}{2})\\&+O(\|\psi_{xxx}\|_{L^\infty(I_1\times I_2)})\Delta^3+O(\|\psi_{xnn}\|_{L^\infty(I_1\times I_2)})\Delta M^2+O(\|\psi_{xxn}\|_{L^\infty(I_1\times I_2)})\Delta^2 M.
\end{align*}
The terms involving $\psi$, $\psi_x$ and $\psi_{xx}$  are independent of $n$, so they are annihilated by the absolute value. The terms involving $\psi_n$ and  $\psi_{nn}$
are independent of $x$. They disappear from the computations once we choose
\begin{equation}
\label{1}
a_n=e(-\psi_n(x_N,\frac{N}2)(n-\frac{N}{2})-\frac12\psi_{nn}(x_N,\frac{N}2)(n-\frac{N}{2})^2 ).\end{equation}
We would like to make the cubic terms $o(1)$. We choose  $\Delta=o(\frac1{N^{\frac{k+1}{3(k-1)}}})$, $M=N^{\frac{2k-1}{3(k-1)}}$. Note that our choice of $x_N$ implies that
$$\psi_{xn}(x_N,\frac{N}2)=0.$$
 Note also that $x\sim x_N$ whenever $x\in I_1$. Thus $$\|\psi_{xnn}\|_{L^\infty(I_1\times I_2)}\lesssim \frac1{N},\;\;\;\text{and }\;\;\|\psi_{xnn}\|_{L^\infty(I_1\times I_2)}\Delta M^2=o(1).$$
Similarly,
$$\|\psi_{xxn}\|_{L^\infty(I_1\times I_2)}\lesssim \frac{N}{N^{\frac{k-2}{k-1}}},\;\;\;\text{and }\;\;\|\psi_{xxn}\|_{L^\infty(I_1\times I_2)}\Delta^2 M=o(1)$$
$$\|\psi_{xxx}\|_{L^\infty(I_1\times I_2)}\lesssim \frac{N^2}{N^{\frac{k-3}{k-1}}},\;\;\;\text{and }\;\;\|\psi_{xxx}\|_{L^\infty(I_1\times I_2)}\Delta^3 =o(1).$$
Let $a_n$ be as in \eqref{1} for $|n-\frac{N}{2}|\le M$, and $a_n=0$ otherwise. So $\|a_n\|_2\sim M^{1/2}$. We have for $|x-x_N|\le \Delta$
$$|\sum_{n=-N}^Na_ne(nx-n^2\phi(x))|\sim M,$$
thus
$$\int_{-\delta}^\delta|\sum_{n=-N}^Na_ne(nx-n^2\phi(x))|^pdx\gtrsim M^p\Delta.$$
It now suffices to note that $M^{1/2}\Delta^{1/p}\sim N^{\frac{2k-1}{6(k-1)}-\frac{k+1}{3p(k-1)}}$.

\end{proof}

\section{$L^2$ upper bounds}
\label{s3}

Fix $k\ge 3$. It suffices to prove that
$$\|\sum_{C_kN\le n\le N}a_ne(nx-n^2\frac{x^k}{k})\|_{L^2(0,1)}\les N^{\frac{k-2}{6(k-1)}}\|a_n\|_2.$$
where $C_k<1$ is close to 1, depending on $k$. Working with $x^k/k$ versus $x^k$ simplifies some of the numerology, but has no other relevance to the argument. The reduction to a fixed dyadic block costs us a $\log N$.  In particular, all future computations will involve $n,m\sim N$. The estimate
$$\|\sum_{-N\le n\le -1}a_ne(nx-n^2\frac{x^k}{k})\|_{L^2(0,1)}\les \|a_n\|_2$$
is an immediate consequence of Lemma \ref{2} and Lemma \ref{3}, as it will become clear in a moment. The advantage here is that there are no critical points in $[0,1]$.
\smallskip

Let $\phi(x)=\phi_{n,m}(x)=(n-m)[x-(n+m)\frac{x^k}{k}]$, with $n\not=m$. This phase function has only one critical point in $[0,1]$, at $x_0=\frac1{(n+m)^{\frac1{k-1}}}$. Let $\theta_{n,m}$ and $\eta$ be smooth functions such that
$$1_{[0,1]}(x)\le \theta_{n,m}(x)+\eta((n+m)^{\frac1{k-1}}(x-x_0))\le 1_{[-x_0/2,2]}(x).$$
In order to facilitate the forthcoming oscillatory integral estimate, we require $\eta$ to be supported in the interval $[-\delta,\delta]$, and to be equal to 1 on $[-\delta/2,\delta/2]$, where $\delta=o(1)$ is small enough depending only on $k$.
\smallskip

We write
\begin{align*}
\int_0^1|\sum_{C_kN\le n\le N}&a_ne(nx-n^2\frac{x^k}{k})|^2dx\le \\&\sum_{C_kN\le n,m\le N}a_n\overline{a_m}\int e((n-m)x-(n^2-m^2)\frac{x^k}{k})\theta_{n,m}(x)dx\\&+\sum_{C_kN\le n,m\le N}a_n\overline{a_m}\int e((n-m)x-(n^2-m^2)\frac{x^k}{k})\eta((n+m)^{\frac1{k-1}}(x-x_0))dx.
\end{align*}

The next two lemmas clarify the asymptotic nature of the integrals in the two sums.

\begin{lem}[Stationary regime]
\label{4}
If $|n-m|\gtrsim N^{\frac1{k-1}}$ we have for some $C_{sta}\in \C$ independent of $n,m$
	$$\int e(\phi(x))\eta((n+m)^{\frac1{k-1}}(x-x_0))dx=O(\frac1{|m-n|})+\begin{cases}\frac{\overline{C_{sta}}}{|n-m|^{\frac12}|n+m|^{\frac1{2(k-1)}}}e(\frac{(k-1)(n-m)}{k(n+m)^{\frac1{k-1}}}), \;\;m<n\\\\\frac{{C_{sta}}}{|n-m|^{\frac12}|n+m|^{\frac1{2(k-1)}}}e(\frac{(k-1)(n-m)}{k(n+m)^{\frac1{k-1}}}), \;\;m>n\end{cases}.$$
	If $|n-m|\lesssim N^{\frac1{k-1}}$ we have$$|\int e(\phi(x))\eta((n+m)^{\frac1{k-1}}(x-x_0))dx|\lesssim \frac1{N^{\frac1{k-1}}}.$$
\end{lem}
\begin{proof}
Use Taylor expansion near $x_0$ to write
$$\phi(x)=\frac{(k-1)(n-m)}{k(n+m)^{\frac1{k-1}}}-\frac{n^2-m^2}{k(n+m)^{\frac{k}{k-1}}}\sum_{j=2}^k{k\choose j}[(n+m)^{\frac1{k-1}}(x-x_0)]^j.$$
Changing variables $(n+m)^{\frac1{k-1}}(x-x_0)=t$ we find that the integral equals
$$e(\frac{(k-1)(n-m)}{k(n+m)^{\frac1{k-1}}})\frac{1}{(n+m)^{\frac1{k-1}}}\int e(\lambda P(t))\eta(t)dt,$$
where $\lambda=\frac{m^2-n^2}{(n+m)^{\frac{k}{k-1}}}$, $P(t)=\frac1k\sum_{j=2}^k{k\choose j}t^j$. When $|\lambda|\lesssim 1$, we use the fact that the integral is $O(1)$.
\smallskip

Let us now assume that $\lambda\gtrsim 1$, so we have $m>n$.
Note that $P(0)=P'(0)=0$ and $P''(0)\not=0$. We may thus invoke Proposition 3 (see also Remarks 1.3.4) in Chapter 8 from \cite{S} to write,  with some $C_{sta}\in \C$ independent of $\lambda$
$$\int e(\lambda P(t))\eta(t)dt=\frac{C_{sta}}{|\lambda|^{1/2}}+O(|\lambda|^{-1}).$$
Note also that when $\lambda$ is negative, $\int e(\lambda P(t))\eta(t)dt$ is the complex conjugate of $\int e(-\lambda P(t))\eta(t)dt$. The proof is now complete.

\end{proof}

\begin{lem}[Non-stationary regime]If $n\not=m$ we have
\label{2}	
$$|\int e(\phi(x))\theta_{n,m}(x)dx|\lesssim \frac1{|n-m|}.$$	

\end{lem}
\begin{proof}
To simplify notation, let us write $\theta=\theta_{n,m}$.	
Since $\phi'$ is nonzero on the support of $\theta$, we can use integration by parts to write
$$|\int e(\phi(x))\theta(x)dx|\le \int|\frac{\theta'(x)}{\phi'(x)}|dx+\int|\frac{\theta(x)\phi''(x)}{(\phi'(x))^2}|dx.$$The support of $\theta$ is inside the union of $I_1=[-x_0/2, x_0(1-\frac\delta2)]$ and $I_2=[x_0(1+\frac\delta2),2]$.
The measure of $I_1$ is $O(x_0)$ and we have $\|\theta'\|_{L^\infty(I_1)}=O(x_0^{-1})$. Also, $|\phi'(x)|\sim |n-m|$ for $x\in I_1$. These show that $\int_{I_1}|\frac{\theta'(x)}{\phi'(x)}|dx=O(\frac1{|n-m|})$.
A similar argument shows that $\int_{I_1}|\frac{\theta(x)\phi''(x)}{(\phi'(x))^2}|dx\lesssim \frac1{|n-m|}$, since  $|\phi''(x)|\lesssim \frac{|n-m|}{x_0}$ on $I_1$.

On the other hand, on $I_2$ we have $\|\theta'\|_{L^\infty(I_2)}=O(x_0^{-1})$, $|\phi'(x)|\sim N|n-m|x^{k-1}$
and $|\phi''(x)|\sim N|n-m|x^{k-2}$. We conclude by writing
$$\int_{I_2}|\frac{\theta'(x)}{\phi'(x)}|dx\lesssim \frac1{N|n-m|x_0}\int_{I_2}\frac1{x^{k-1}}dx\lesssim \frac1{|n-m|},$$
$$\int_{I_2}|\frac{\theta(x)\phi''(x)}{(\phi'(x))^2}|dx\lesssim \frac1{N|n-m|}\int_{I_2}x^{-k}dx\lesssim \frac1{|n-m|}.$$

\end{proof}
The bound
$$\int|\sum_{C_kN\le n\le N}a_ne(nx-n^2\frac{x^k}{k})|^2\theta(x)dx\lesssim (\log N)\|a_n\|_2^2$$
is a consequence of Lemma \ref{2} and Lemma \ref{3}.
\smallskip

It remains to prove the bound
$$\int|\sum_{CN\le n\le N}a_ne(nx-n^2\frac{x^k}{k})|^2\eta((n+m)^{\frac1{k-1}}(x-x_0))dx\les N^{\frac{k-2}{3(k-1)}} \|a_n\|_2^2. $$
By invoking Lemma \ref{4} and Lemma \ref{3} (to deal with the error term $O(\frac1{|m-n|})$), it suffices to prove that
\begin{equation}
\label{5}
|\sum_{C_kN\le n,m\le N}a_n\overline{a_m}c_{n,m}|\les N^{\frac{k-2}{3(k-1)}}\|a_n\|_2^2,\end{equation}
where
$$c_{n,m}=\begin{cases}\frac{\overline{C_{sta}}}{|n-m|^{\frac12}|n+m|^{\frac1{2(k-1)}}}\;e(\frac{(k-1)(n-m)}{k(n+m)^{\frac1{k-1}}}),\;\text{if}\;|n-m|\ge N^{\frac1{k-1}}\text{ and }m<n\\\\\frac{{C_{sta}}}{|n-m|^{\frac12}|n+m|^{\frac1{2(k-1)}}}\;e(\frac{(k-1)(n-m)}{k(n+m)^{\frac1{k-1}}}),\;\text{if}\;|n-m|\ge N^{\frac1{k-1}}\text{ and }m>n\\\\0,\;\text{otherwise}\end{cases}.$$
We first note that Lemma \ref{3} is not strong enough to give \eqref{5}, as for each $n\sim N$ we have
$$\sum_{m:\;C_kN\le m\le N\atop{|n-m|\ge N^{\frac1{k-1}}}}\frac{1}{|n-m|^{\frac12}|n+m|^{\frac1{2(k-1)}}}\sim N^{\frac{k-2}{2(k-1)}}.$$
To get the stronger result we need to exploit the oscillatory nature of $c_{n,m}$, and Lemma \ref{6} will prove the right tool for this purpose. It is immediate that $C=(c_{n,m})$
is Hermitian. We need to estimate the coefficients
$$d_{n,m}=\sum_{x}c_{n,x}{c_{x,m}}.$$

\smallskip

Given $C_kN\le n< m\le N$, consider the function $f=f_{n,m}$ defined for $C_kN\le x\le N$
$$f(x)=\frac{k-1}{k}(\frac{x-n}{(x+n)^{\frac1{k-1}}}-\frac{x-m}{(x+m)^{\frac1{k-1}}}).$$
We start by proving two auxiliary lemmas regarding the derivatives of $f$. The coefficient $\frac{k-1}k$ will be irrelevant for our computations.

\begin{lem}
\label{7}	
We have for each $C_kN\le x\le N$ $$|f''(x)|\sim \frac{|m-n|}{N^{\frac{2k-1}{k-1}}}.$$	
\end{lem}
\begin{proof}
We compute
\begin{equation}
\label{9}
f'(x)=\frac1{k}\left[\frac{(k-2)x+kn}{(x+n)^{\frac{k}{k-1}}}-\frac{(k-2)x+km}{(x+m)^{\frac{k}{k-1}}}\right],
\end{equation}
$$f''(x)=\frac{1}{k(k-1)}\left[\frac{(2-k)x+(2-3k)n}{(x+n)^{\frac{2k-1}{k-1}}}-\frac{(2-k)x+(2-3k)m}{(x+m)^{\frac{2k-1}{k-1}}}\right].$$
Define $$g_x(t)=\frac{(2-k)x+(2-3k)t}{(x+t)^{\frac{2k-1}{k-1}}}$$ and note that
$$g_x'(t)=\frac{-k^2x+tk(3k-2)}{(k-1)(x+t)^{\frac{3k-2}{k-1}}}.$$
Note that $k^2<k(3k-2)$ when $k\ge 2$. Choosing $C_k$ sufficiently close to 1 insures that for $C_kN\le x,t\le N$
we have
$$g_x'(t)\sim N^{-\frac{2k-1}{k-1}}.$$
The desired estimate now follows from the Mean Value Theorem.
	
\end{proof}	
\begin{lem}
\label{8}
If $|x-n|\ll m-n$ or $|x-m|\ll m-n$ then
\begin{equation}
\label{14}
|f'(x)|\sim \frac{(m-n)^2}{N^{\frac{2k-1}{k-1}}}.\end{equation}
\end{lem}
\begin{proof}
We start by pointing out that the restriction on $x$ being close to either $n$ or $m$ is necessary. Indeed, since $f(n)=f(m)$, $f$ must have a critical point in the interval $(n,m)$. The lemma says that this point is somewhere in the middle. 	
	
Let $$h_x(t)=\frac{(k-2)x+kt}{(x+t)^{\frac{k}{k-1}}}.$$
Then
$$h_x'(t)=\frac{k}{k-1}\times \frac{x-t}{(x+t)^{\frac{2k-1}{k-1}}}.$$
Since $h_x'$ can be very small (for example $h_x'(x)=0$), the Mean Value Theorem is not useful in estimating \eqref{9}. Instead, we first use Taylor's formula with cubic error terms to  estimate $|f'(n)|$ and $|f'(m)|$. More precisely,
\begin{equation}
\label{10}
f'(m)=h_m(n)-h_m(m)=\frac12h_m''(m)(n-m)^2+O(\|h_m'''\|_{L^\infty [C_kN,N]})(m-n)^3.
\end{equation}
We compute
$$h_x''(t)=\frac{k}{(k-1)^2}\times \frac{-kx+(3k-2)t}{(x+t)^{\frac{3k-2}{k-1}}},$$
$$h_x'''(t)=\frac{k(2k-1)(3k-2)}{(k-1)^3}\times \frac{x-t}{(x+t)^{\frac{4k-3}{k-1}}}.$$
The relevant things for us are $h_m''(m)\sim \frac1{N^{\frac{2k-1}{k-1}}}$ and $\|h_m'''\|_{L^\infty [C_kN,N]}\lesssim \frac1{N^{\frac{3k-2}{k-1}}}$. By choosing $C_k$ sufficiently close to 1 we can guarantee that $m-n=o(N).$ This in turn makes the cubic error term in \eqref{10} much smaller than the quadratic term, showing that
$$|f'(m)|\sim \frac{(m-n)^2}{N^{\frac{2k-1}{k-1}}}.$$
Finally, using Lemma \ref{7} and the Mean Value Theorem we find that
$$|f'(x)-f'(m)|\sim |x-m|\frac{m-n}{N^{\frac{2k-1}{k-1}}}.$$
This quantity is $o(\frac{(m-n)^2}{N^{\frac{2k-1}{k-1}}})$ when $|x-m|\ll m-n$, showing that $|f'(x)|\sim \frac{(m-n)^2}{N^{\frac{2k-1}{k-1}}}.$ A similar argument works when $x$ is near $n$.

\end{proof}

We are now ready to estimate the coefficients $d_{n,m}$. We start with the easy regime where there are no cancellations.
\begin{pr}
\label{12}	
Assume $|n-m|\le N^{\frac{2k-1}{3(k-1)}}$. Then $|d_{n,m}|\les \frac1{N^{\frac1{k-1}}}$.
\end{pr}
\begin{proof}Triangle inequality leads to the estimate
$$|d_{n,m}|\lesssim \frac1{N^{\frac1{k-1}}}\sum_{x\sim N\atop{x\not=n,m}}\frac{1}{\sqrt{|n-x|}}\frac1{\sqrt{|m-x|}}.$$
When $|x-\frac{n+m}{2}|\lesssim |m-n|$, we use that $\max\{|x-n|,|x-m|\}\gtrsim |n-m|$. When $|x-\frac{n+m}{2}|\gg |m-n|$ we use that $|x-n|,|x-m|\sim |x-\frac{n+m}{2}|$.

\end{proof}

Before we move to discuss  the next regime, let us pause to understand the sharpness in the estimate from Proposition \ref{12}. Let us assume $|n-m|\sim N^{\frac{2k-1}{3(k-1)}}$. Let $x_0\in (n,m)$ be the critical point of $f=f_{n,m}$. Then, using Lemma \ref{7} we write for each $|x-x_0|\ll m-n$
$$f(x)=f(x_0)+O(\frac{m-n}{N^{\frac{2k-1}{k-1}}})(x-x_0)^2=f(x_0)+o(1).$$
It follows that when evaluating the sum coming from such $x$ (as part of the whole sum defining $d_{n,m}$), there is no cancellation. More precisely
$$|\sum_{|x-x_0|\ll m-n}\frac{1}{\sqrt{|(n-x)(m-x)|}}\frac1{[(n+x)(m+x)]^{\frac1{2(k-1)}}}e(f(x))|$$$$\sim \sum_{|x-x_0|\ll m-n}\frac{1}{\sqrt{|(n-x)(m-x)|}}\frac1{[(n+x)(m+x)]^{\frac1{2(k-1)}}}\sim \frac1{N^{\frac1{k-1}}}.$$
\smallskip

We will next see that cancellations occur when $m$ and $n$ are further apart.
\begin{pr}
\label{17}	
If $|n-m|\ge N^{\frac{2k-1}{3(k-1)}}$ we have
$|d_{n,m}|\lesssim \frac{N^{\frac{2k-3}{2(k-1)}}}{|n-m|^{3/2}}$.
\end{pr}
\begin{proof}
We may assume that $n<m$.	
We distinguish four regimes. Let $I_1$, $I_2$ be  intervals  of length $o(m-n)$ centered at $n$ and $m$ respectively, such that \eqref{14} holds. Let $J_1$, $J_2$ be  intervals  of smaller length $o( \frac{N^{\frac{2k-1}{k-1}}}{|n-m|^2})$ centered at $n$ and $m$ respectively. Note that $J_1\subset I_1$ and $J_2\subset I_2$.

The choice of $J_1$ and $J_2$ is motivated by the fact that there is no cancellation for the sums of  $e(f(x))$ when $x$ is restricted to these intervals. Indeed, Lemma \ref{8} shows that when $x\in J_1$
$$|f(x)-f(n)|\sim \frac{|J_1|(m-n)^2}{N^{\frac{2k-1}{k-1}}}\ll 1.$$

As a general rule throughout the forthcoming argument, we will use the fact that the sequence  $$w_x=\frac{1}{\sqrt{|(n-x)(m-x)|}}\frac1{[(n+x)(m+x)]^{\frac1{2(k-1)}}}$$
changes monotonicity $O(1)$ many times in the interval $[C_kN,N]$. Thus, for each interval $I$ we have that $\|w_x\|_{V^1(I)}\lesssim \|w_x\|_{l^\infty(I)}$.

The reader will notice that the difference between $C_{sta}$ and $\overline{C_{sta}}$ in the definition on $c_{n,m}$ does not affect the argument.
\\
\\
1.  Let $x$ be in the core of $[n,m]$, more precisely $x\in R_1=[n,m]\setminus (I_1\cup I_2)$. This interval contains the critical point of $f$, so the first derivative test is not applicable. However, we may use the second derivative test, Theorem \ref{15} to estimate for each $I\subset R_1$
$$|\sum_{x\in I}e(f(x))|\lesssim \lambda^{1/2}|I|+\lambda^{-1/2},$$
where, according to Lemma \ref{7} we have $\lambda\sim \frac{m-n}{N^{\frac{2k-1}{k-1}}}$. The term $\lambda^{-1/2}$ dominates since $|I|\le m-n\lesssim \lambda^{-1}$. Moreover,  $\|w_x\|_{l^{\infty}(R_1)}\sim \frac1{N^{\frac1{k-1}}(m-n)}$. Combining this with Lemma \ref{16} we are led to the desired estimate
$$|\sum_{x\in R_1}w_xe(f(x))|\lesssim \lambda^{-1/2}\|w_x\|_{l^{\infty}(R_1)}\sim \frac{N^{\frac{2k-3}{2(k-1)}}}{|n-m|^{3/2}}.$$
\\
\\
2. An identical argument works for $R_2=[C_kN,N]\setminus (I_1\cup I_2\cup [n,m])$. The term $\lambda^{-1/2}$ still dominates and we have the same upper bound for $\|w_x\|_{l^{\infty}(R_2)}$ as before. While we could do better in this regime by using the first derivative test, the previous upper bound is good enough for our purposes.
\\
\\
3. Let us now restrict $x$ to one of the four intervals comprising the set $R_3=(I_1\setminus J_1)\cup(I_2\setminus J_2)$. Note first that the $l^\infty$ norm of the weights gets larger in this regime, due to the proximity of $x$ to either $n$ or $m$
$$\|w_x\|_{l^\infty(R_3)}\sim \frac{1}{N^{\frac1{k-1}}\sqrt{m-n}\sqrt{|J_1|}}\sim \frac{\sqrt{m-n}}{N^{\frac{2k+1}{2(k-1)}}}.$$
An application of the second derivative test would lead to the bound
$$|\sum_{x\in R_3}w_xe(f(x))|\lesssim \frac{\sqrt{m-n}}{N^{\frac{2k+1}{2(k-1)}}}\lambda^{-1/2}\sim \frac1{N^{\frac1{k-1}}}.$$
However,  this is significantly worse than what we need to prove. To achieve the desired result, we use the first derivative test
$$|\sum_{x\in R_3}w_xe(f(x))|\lesssim \frac{\sqrt{m-n}}{N^{\frac{2k+1}{2(k-1)}}}\times \frac{N^{\frac{2k-1}{k-1}}}{(m-n)^2}\sim \frac{N^{\frac{2k-3}{2(k-1)}}}{(m-n)^{3/2}}. $$
\\
\\
4. When $x\in R_4=J_1\cup J_2$, the sharp bound follows using the triangle inequality
$$|\sum_{x\in R_4}w_xe(f(x))|\lesssim \sum_{x\in R_4}w_x\sim \frac1{N^{\frac1{k-1}}\sqrt{m-n}}\sum_{1\le p\ll \frac{N^{\frac{2k-1}{k-1}}}{|n-m|^2}}\frac1{\sqrt{p}}\sim \frac{N^{\frac{2k-3}{2(k-1)}}}{(m-n)^{3/2}}.$$

\end{proof}

We are finally in position to prove \eqref{5}. For each $n\in [C_kN,N]$, we use Lemma \ref{12} and \ref{17} to write
\begin{align*}
\sum_{C_kN\le m\le N}|d_{n,m}|&\le \sum_{m:\;|m-n|\le N^{\frac{2k-1}{3(k-1)}}}|d_{n,m}|+ \sum_{m:\;|m-n|\ge N^{\frac{2k-1}{3(k-1)}}}|d_{n,m}|\\&\les N^{\frac{2k-1}{3(k-1)}-\frac1{k-1}}+{N^{\frac{2k-3}{2(k-1)}}}\sum_{N^{\frac{2k-1}{3(k-1)}}\le p\lesssim N}\frac1{p^{3/2}}\\&\sim N^{\frac{2(k-2)}{3(k-1)}}.\end{align*}

We combine this estimate with Lemma \ref{6} to write
\begin{align*}
\sum_{C_kN\le n,m\le N}a_n\overline{a_m}c_{n,m}\le (\max_n\sum_m|d_{n,m}|)^{1/2}\|a_n\|_2^2\lesssim N^{\frac{k-2}{3(k-1)}}\|a_n\|_2^2.
\end{align*}
\smallskip

\begin{re}
It seems possible that the upper bound (in addition to the already mentioned lower bound) in Theorem \ref{20} would hold for all $C^{k+1}$ functions with finitely many zeros of $\phi'$, all of order  $k-1$. To keep things concise, we have decided to not pursue this line of investigation. However, there is a simple argument for $k=3$, that we mention briefly. Let us assume $\phi'(0)=\phi''(0)=0$ and $\phi'''(0)\not=0$, so that $\phi(x)=\frac{\phi'''(0)}{3!}x^3+O(x^4)$. We may assume 0 is the only critical point of $\phi$.

The estimate in the nonstationary regime
$$\|\sum_{n\sim N}a_ne(nx+n^2\phi(x))\|_{L^2(|x|\gg N^{-1/2})}\les \|a_n\|_2$$
is an immediate consequence of the fact that the critical point of  $$\phi_{n,m}(x)=(n-m)[x+(n+m)\phi(x)]$$ lies in the interval $|x|\lesssim N^{-1/2}$ when $n\not=m\sim N$. Note further that we have 
\begin{equation}
\label{jfiurng9iug80uy0}
\sum_{n=1}^N\|w_{n+1}(x)-w_n(x)\|_{L^\infty(|x|\lesssim N^{-1/2})}\lesssim 1,
\end{equation}
 where $w_n(x)=e(n^2(\phi(x)-\frac{\phi'''(0)}{3!}x^3))$. Thus, using summation by parts and  Cauchy--Schwarz we have
$$\|\sum_{n=1}^Na_ne(nx+n^2\phi(x))\|_{L^2(|x|\lesssim N^{-1/2})}^2\lesssim$$
$$\int_{|x|\lesssim N^{-1/2}}(\sum_{m=1}^N|w_{m+1}(x)-w_m(x)|)(\sum_{m=1}^N|w_{m+1}(x)-w_m(x)||\sum_{n=1}^ma_ne(nx+n^2\frac{\phi'''(0)}{6}x^3)|^2)dx\lesssim$$
$$\|\sum_{m=1}^N|w_{m+1}(x)-w_m(x)|\|_{L^\infty(|x|\lesssim \frac1{N^{1/2}})}\int_{|x|\lesssim \frac1{N^{1/2}}}\sum_{m=1}^N|w_{m+1}(x)-w_m(x)||\sum_{n=1}^ma_ne(nx+n^2\frac{\phi'''(0)}{6}x^3)|^2dx\lesssim$$

$$(\sum_{m=1}^N\|w_{m+1}(x)-w_m(x)\|_{L^\infty(|x|\lesssim N^{-1/2})})^2\max_{m\le N}\int_{|x|\lesssim N^{-1/2}}|\sum_{n=1}^ma_ne(nx+n^2\frac{\phi'''(0)}{6}x^3)|^2dx.$$
This is $\les N^{1/6}\|a_n\|_2^2$, due to Theorem \ref{20} and \eqref{jfiurng9iug80uy0}.

\end{re}

\section{Related results and open questions}
\label{s4}

The two dimensional analogue of Theorem \ref{23} proved in \cite{DL} is as follows.
\begin{te}
\label{24}
Assume $\phi:[-1,1]^2\to\R$ is $C^2$ and has nonzero Hessian $$\inf_{(x,y)\in[-1,1]^2}|\phi_{xx}(x,y)\phi_{yy}(x,y)-\phi_{xy}^2(x,y)|>0.$$Then for each $a_{n,m}\in\C$ we have
$$\|\sum_{n,m=-N}^Na_{n,m}e(nx+my+(n^2+m^2)\phi(x,y))\|_{L^2([-1,1]^2)}\les \|a_{n,m}\|_{2}.$$	
\end{te}
Let us show that the result may fail for the zero curvature cylinder $\phi(x,y)=-x^3$.
\begin{te}
For each $N$ there is a choice of coefficients $a_{n,m}$ such that
$$\|\sum_{n,m=-N}^Na_{n,m}e(nx+my-(n^2+m^2)x^3)\|_{L^2([-1,1]^2)}\gtrsim N^{\frac1{12}} \|a_{n,m}\|_{2}.$$	
\end{te}
\begin{proof}
Write the phase as
$$(nx-n^2x^3)+(mx-m^2x^3)+(m-\frac{N}{2})(y-x)+\frac{N(y-x)}{2}.$$
The proof of Theorem \ref{21} shows that for appropriate functions $f_1,f_2$ and for $|n-\frac{N}{2}|,|m-\frac{N}{2}|\ll N^{5/6}$, $|x-\frac1{\sqrt{3N}}|\ll \frac1{N^{2/3}}$ we have
$$nx-n^2x^3=f_1(n)+f_2(x)+o(1),$$
$$mx-m^2x^3=f_1(m)+f_2(x)+o(1).$$
Let $a_{n,m}=e(-f_1(n)-f_1(m))$ if $|n-\frac{N}{2}|,|m-\frac{N}{2}|\ll N^{5/6}$, and $a_{n,m}=0$ otherwise.
It follows that we have constructive interference
$$|\sum_{n,m=-N}^Na_{n,m}e(nx+my-(n^2+m^2)x^3)|\sim N^{5/3}$$
on the set $ S=\{(x,y):\;|x-\frac1{\sqrt{3N}}|\ll \frac1{N^{2/3}},\;|y-x|\ll N^{-5/6}\}$. The result is now immediate.

\end{proof}	
The proof in \cite{DL} of Theorem \ref{24}  uses solely the decay of the Fourier transform of the surface measure of the graph of $\phi$, but not its oscillatory component. This leads to a logarithmic loss in the upper bound, and to even larger losses in higher dimensions. It seems possible that a version of the $TT^*$ approach developed here could remove these losses.
\smallskip

On a related note, we mention a result of Bourgain and Rudnick for exponential sums supported on dilates of the sphere$$S^{\S}_{a,d}(\x,N)=\sum_{\n\in\sqrt{N}\S^{d-1}\cap \Z^d}a_{\n}e(\n\cdot \x).$$
\begin{te}[\cite{BR}]
Let $\Sigma$ be a real analytic hypersurface in $\T^d$, with nonzero curvature and surface measure $\sigma$. Then if $d=2$ or $d=3$ we have
	\begin{equation}
	\label{a40}
	\|S^{\S}_{a,d}(\x,N)\|_{L^2(d\sigma)}\lesssim \|a\|_{2}.
	\end{equation}
\end{te}
In light of our examples for the parabola and the paraboloid, the following question is natural. See also Conjecture 1.9 in \cite{BR}.
\begin{qu}Is there a real analytic $\Sigma$ such that \eqref{a40} is false?
\end{qu}	
\smallskip

The result for the parabola $(n,n^2)$ in Theorem \ref{23} admits yet another extension to three dimensions. Rather than considering the paraboloid as in Theorem \ref{24}, one may consider the moment curve $(n,n^2,n^3)$. The following result was proved in \cite{DL}.
\begin{te}
	\label{26}
	Assume $\phi:[-1,1]^2\to\R$  has nonzero Hessian. Then for each $a_{n}\in\C$ we have
	\begin{equation}
	\label{60}
	\|\sum_{n=-N}^Na_{n}e(nx+n^2y+n^3\phi(x,y))\|_{L^6([-1,1]^2)}\lesssim_\epsilon N^\epsilon \|a_{n}\|_{2}.
	\end{equation}	
\end{te}
The example $\phi(x,y)=xy$ produces constructive interference near the origin ($a_n=1$), showing that the exponent 6 cannot be improved.
The question left open in \cite{DL} is about whether the above inequality could hold for arbitrary $C^\infty$ functions $\phi$. It was observed that the result does indeed hold when $\phi=0$. This may suggest that perhaps the properties of $\phi$ should play no serious role in the argument. The following result proves however that curvature is critical.
\begin{te}
\label{31}	
For each $N$ there is a sequence $a_n$ such that
$$\|\sum_{n=-N}^Na_ne(nx+n^2y-n^3(y+x^3)^3)\|_{L^6([-1,1]^2)}\gtrsim N^{1/36}\|a_n\|_2.$$
\end{te}
\begin{proof}
With the change of variables $u=x$, $v=y+x^3$, we have
$$\|\sum_{n=-N}^Na_ne(nx+n^2y-n^3(y+x^3)^3)\|_{L^6([-1,1]^2)}\ge$$$$\|\sum_{n=-N}^Na_ne((nu-n^2u^3)+(n^2v-n^3v^3))\|_{L^6([-1/2,1/2]^2)}.$$
Recall that if $|n-\frac{N}{2}|\ll N^{5/6}$ and $|u-\frac1{\sqrt{3N}}|\ll N^{-2/3}$ we have
$$nu-n^2u^3=f_1(n)+f_2(u)+o(1).$$
A similar argument will show that with $v_N=\frac{2}{3\sqrt{N}}$ we have
$$\theta(v,n)=n^2v-n^3v^3+n^3v_N^3=g_1(n)+g_2(v)+o(1)$$
whenever  $(v,n)\in J_1\times J_2$ defined by $|v-v_N|\ll \delta=N^{-1-\frac23}$ and $|n-\frac{N}{2}|\ll M=N^{5/6}$. This follows as before from Taylor's formula with cubic error terms, since
$$\theta_{vn}(v_N,\frac{N}2)=0,$$
 $$\|\theta_{vnn}\|_{L^\infty(J_1\times J_2)}\lesssim 1,\;\;\;\text{and }\;\;\|\theta_{vnn}\|_{L^\infty(J_1\times J_2)}\delta M^2=o(1),$$
$$\|\theta_{vvn}\|_{L^\infty(J_1\times J_2)}\lesssim N^{3/2},\;\;\;\text{and }\;\;\|\theta_{vvn}\|_{L^\infty(J_1\times J_2)}\delta^2 M=o(1),$$
$$\|\theta_{vvv}\|_{L^\infty(J_1\times J_2)}\lesssim N^3,\;\;\;\text{and }\;\;\|\theta_{vvv}\|_{L^\infty(J_1\times J_2)}\delta^3 =o(1),$$
and
$$\|\theta_{nnn}\|_{L^\infty(J_1\times J_2)}\lesssim \max_{v\in J_1}|v^3-v_N^3|\sim \frac{1}{N^{8/3}},\;\;\;\text{and }\;\;\|\theta_{nnn}\|_{L^\infty(J_1\times J_2)}M^3 =o(1).$$
It follows that for an appropriate choice of $a_n$ restricted to  $|n-\frac{N}{2}|\ll N^{5/6}$ we have
$$\int_{|u-\frac1{\sqrt{3N}}|\ll N^{-2/3}\atop{|v-\frac{2}{3\sqrt{N}}|\ll N^{-5/3}}}|\sum_{n=-N}^Na_ne((nu-n^2u^3)+(n^2v-n^3v^3))|^6dudv\sim \frac{N^5}{N^{7/3}}\sim N^{1/6}\|a_n\|_2^6.$$
\end{proof}
Another interesting question is whether \eqref{60} would hold for some $L^p$, $p>6$, for some, or maybe all surfaces with both principal curvatures positive. For example, constructive interference similar to the one discussed here does not seem to rule out the possibility that the following inequality could be true
$$\|\sum_{n=-N}^Na_n(nx+n^2y+n^3(x^2+y^2))\|_{L^{7}([-1,1]^2)}\lesssim_\epsilon N^\epsilon \|a_n\|_{2}.$$
Even the case $a_n=1$ seems interesting.
If this inequality does indeed hold, it seems that the methods developed in \cite{DL} and here would not be enough to prove it. Quite likely,  the proof would need to rely on estimates for cubic Weyl sums, not just on  Fourier analysis. This probably remains true even for the following special case obtained by using $|x|\ll N^{-3/2}$
$$\|\sum_{n=-N}^Na_n(n^2y+n^3y^2)\|_{L^{7}([-1,1]^2)}\lesssim_\epsilon N^{\frac3{14}+\epsilon} \|a_n\|_{2}.$$ 
\smallskip

Theorem \ref{26} provides a positive answer in three dimensions to the following conjecture proposed in \cite{DL}.
\begin{con}
\label{30}	
If $\Sigma$ is a smooth hypersurface in $\R^d$ with nonzero Gaussian curvature and surface measure $d\sigma$, then	for $p\le p_{d-1}=d(d-1)$
\begin{equation}
\label{29}
\|\sum_{n=-N}^Na_ne((n,n^2,\ldots,n^d)\cdot \x)\|_{L^{p}(d\sigma)}\lesssim_\epsilon N^\epsilon\|a_n\|_2.
\end{equation}	
\end{con}
While partial results have been achieved in \cite{DL}, the conjecture has only been verified in the full range when $d=3$. Interestingly, inequality \eqref{29} is true if $\Sigma$ is the (zero curvature) hyperplane $x_d=0$. This is the weighted version of the Main Conjecture in Vinogradov's Mean Value Theorem first proved in \cite{BDG}, and then reproved in several other ways, \cite{Wol},\cite{GLYZ}. It seems reasonable to ask whether one of these arguments would apply to prove Conjecture \ref{30}. Typically, these arguments use induction on dimension. To prove \eqref{29} in dimension $d$, one may need a similar essentially scale-independent bound in $d-1$ dimensions, for $p=p_{d-2}$ and $\Sigma'$ a $d-2$ dimensional slice of $\Sigma$ (intersection of $\Sigma$ with a hyperplane). If $\Sigma$ is only assumed to have nonzero curvature, its slices might have arbitrarily small curvature. Examples such as the one described in Theorem \ref{21} or Theorem \ref{31} exist in all dimensions, and show the failure of square root cancellation in the absence of curvature. This seems to indicate that an argument along these lines is not feasible for attacking  Conjecture \ref{30}. There is a possibility however that such an argument could be carried out for the class of hypersurfaces with positive principal curvatures, as their slices inherit this property.

\smallskip

We close with a few remarks on other $L^p$ spaces.

\begin{te}
	\label{44}
Let $k\ge 2$ and $p\ge 2$. Then	
$$\|\sum_{n=-N}^Na_ne(nx-n^2x^k)\|_{L^p(-1,1)}\lesssim_\epsilon \|a_n\|_{2}\begin{cases}N^{\frac{2k-1}{6(k-1)}-\frac{k+1}{3p(k-1)}+\epsilon},\;2\le p\le 4\\N^{\frac12-\frac1p+\epsilon},\;p\ge 4\end{cases}.$$
Apart from $N^\epsilon$, the upper bound is sharp.
\end{te}
\begin{proof}
When $k=2$, the upper bound follows by interpolating Theorem \ref{23} with the trivial $L^\infty$ bound. The lower bound is seen via constructive interference on the set $|x|\le \frac1N$, with $a_n=1$.	
	
Let us now focus on $k\ge 3$.	
The lower bound is provided by Theorem \ref{21} when $p\in [2,4]$ and by using constructive interference on the set $|x|\le \frac1N$  when $p\ge 4$.	
	
Due to Theorem \ref{20} and the trivial $L^\infty$ bound, it suffices to prove the $L^4$ estimate
$$\|\sum_{n=-N}^Na_ne(nx-n^2x^k)\|_{L^4(-1,1)}\lesssim_\epsilon N^{\frac14+\epsilon}\|a_n\|_{2}.$$
The argument in Section \ref{s3} shows that
$$\|\sum_{n=-N}^Na_ne(nx-n^2x^k)\|_{L^2(|x|\gg N^{-\frac{1}{k-1}})}\les \|a_n\|_{2}.$$
Thus, the desired $L^4$ estimate in this range follows by interpolation (H\"older)  with the $L^\infty$ bound. On the other hand, since $|x^k|\lesssim 1/N$ in the remaining range $|x|\lesssim N^{-\frac{1}{k-1}}$, it suffices to prove the stronger inequality
$$\|\sup_{|t|\lesssim 1/N}|\sum_{n=-N}^Na_ne(nx+n^2t)|\|_{L^4(-1,1)}\lesssim_\epsilon N^{\frac14+\epsilon}\|a_n\|_{2},$$
or equivalently,
$$\|\sup_{|t|\lesssim N}|\sum_{n=-N}^Na_ne(\frac{n}{N}x+(\frac{n}{N})^2t)|\|_{L^4(-N,N)}\lesssim_\epsilon N^{\frac12+\epsilon}\|a_n\|_{2}.$$
This however is a particular instance of the estimate from \cite{KPV}
$$\|\sup_{|t|\le N}|Ef(x,t)|\|_{L^4(-N,N)}\lesssim_\epsilon N^\epsilon\|f\|_{L^2(-1,1)}$$
for the extension operator $$Ef(x,t)=\int f(\xi)e(\xi x+\xi^2 t)d\xi,$$
using
$$f(\xi)=\sum_{n=-N}^Na_n1_{[\frac{n}{N},\frac{n+1}{N}]}(\xi).$$
 \end{proof}

It is an easy application of nonstationary phase analysis (or alternatively, change variables and use Proposition \ref{oiuur8ur89u98}) to show that if $\inf_{|x|\le 1}|\phi'(x)|>0$ then
$$\|\sum_{n=-N}^Na_ne(nx+n^2\phi(x))\|_{L^2(-1,1)}\les \|a_n\|_{2}.$$
The following seems reasonable to ask, and is known to have a positive answer when $\phi$ is linear.
\begin{qu}
Is it true that if $\inf_{|x|\le 1}|\phi'(x)|>0$ then
$$\|\sum_{n=-N}^Na_ne(nx+n^2\phi(x))\|_{L^4(-1,1)}\lesssim_\epsilon N^\epsilon \|a_n\|_{2}?$$	
\end{qu}
In fact, even the following stronger version of this has a good  chance of being true.
\begin{qu}
	Is it true that given any measurable $\phi:[0,1]\to[0,1]$ we have
	$$\|\sum_{n=-N}^Na_ne(n\phi(x)+n^2x)\|_{L^4(-1,1)}\lesssim_\epsilon N^\epsilon \|a_n\|_{2}?$$	
\end{qu}
Even the $L^2$ estimate seems difficult, although we can easily prove the following.
\begin{pr}
	\label{oiuur8ur89u98}
If $\phi$ is Lipschitz then
$$\|\sum_{n=N/2}^Na_ne(n\phi(x)+n^2x)\|_{L^2(-1,1)}\lesssim  \|a_n\|_{2}.$$
\end{pr}
\begin{proof}
Since $n(\phi(x)-\phi(y))=O(1)$ when $|x-y|\lesssim 1/N$, we may assume $\phi$ is constant at this scale.  It thus suffices to prove that for each interval $I$ with length $1/N$ we have
$$\|\sum_{n=N/2}^Na_ne(n^2x)\|_{L^2(I)}\lesssim  |I|^{1/2}\|a_n\|_{2},$$
or equivalently, that 
$$\|\sum_{n=N/2}^Na_ne((\frac{n}{N})^2x)\|_{L^2(N^2I)}\lesssim N^{1/2}\|a_n\|_{2}.$$
This is a consequence of $L^2$ orthogonality, since the frequencies $(\frac{n}{N})^2$ are $1/N$-separated and $N^2I$ has length $N$.

\end{proof}

\section{Auxiliary results}
The proofs of the following derivative tests can be found in \cite{GK}.

\begin{te}[First derivative test]
	Assume $I\subset \R$ is an interval of length at least one. Assume $f\in C^1(I)$, $f'$ is monotonic on $I$ and $|f'(x)|\sim \lambda$ for some $\lambda=o(1)$. Then
	$$|\sum_{x\in I\cap\Z}e(f(x))|\lesssim \lambda^{-1}.$$	
\end{te}

\begin{te}[Second derivative test]
	\label{15}
	Assume $I\subset \R$ is an interval of length at least one. Assume $f\in C^2(I)$,  and $|f''(x)|\sim \lambda$ for some $\lambda>0$. Then
	$$|\sum_{x\in I\cap\Z}e(f(x))|\lesssim \lambda^{1/2}|I|+\lambda^{-1/2}.$$	
\end{te}

Let $C=(c_{n,m})$ be a square matrix with complex entries. We will use the following well known upper bounds for its  spectral norm.
\begin{lem}
	\label{3}
	For each $a_n\in \C$
	$$|\sum_{n,m}a_n\overline{a_m}c_{n,m}|\le \|a_n\|_2^2(\max_n\sum_m|c_{n,m}|+\max_m\sum_n|c_{n,m}|).$$	
\end{lem}	

\begin{lem}
	\label{6}
	Assume $C$ is Hermitian, that is $c_{n,m}=\overline{c_{m,n}}$. Let $$d_{n,m}=\sum_{x}c_{n,x}{c_{x,m}}.$$
	Then for each $a_n\in \C$
	$$|\sum_{n,m}a_n\overline{a_m}c_{n,m}|\le \|a_n\|_2^2(\max_n\sum_m|d_{n,m}|)^{1/2}.$$	
\end{lem}

We control weighted exponential sums over intervals using summation by parts.
\begin{lem}
	\label{16}
	Let $I=\{a+1,\ldots,b\}$ be an interval in $\Z$ and let $f:I\to\R.$ Then for each $w_x\in\R$
	$$|\sum_{n\in I}w_xe(f(x))|\le \|w_x\|_{V^1(I)}\sup_{a+1\le n\le b}|\sum_{a+1\le x\le n}e(f(x))|.$$
\end{lem}

\end{document}